\documentclass{article}

\title{A generalization of Clausen's identity}

\author{Raimundas Vid\=unas\\
 \em Kobe University}

\newtheorem{theorem}{Theorem}[section]

\newtheorem{conjecture}[theorem]{Conjecture}
\newtheorem{lemma}[theorem]{Lemma}

\usepackage{latexsym}

\newcommand{\ediff}[1]{\Theta_#1}
\newcommand{\comm}{\!\cdot\!}
\newcommand{\ahalf}{{\textstyle\frac12}}
\newcommand{\app}[4]{F_{\!#1}\!
  \left(\left.{#2 \atop #3}\right| #4 \right) }

\newcommand{\hpg}[5]{{}_{#1}\mbox{\rm F}_{\!#2}\!
  \left(\left.{#3 \atop #4}\right| #5 \right) }

\newcommand{\hpgo}[2]{{}_{#1}\mbox{\rm F}_{\!#2}}

\newcommand{\proof}{{\bf Proof. }}
\newcommand{\qed}{\hfill $\Box$}
\newcommand{\equal}{&\!\!\!=\!\!\!&}
\newcommand{\rhallign}{\hspace{-26pt}&&} 
\newcommand{\rhallignh}{\rhallign\hspace{50pt}} 

\newcommand{\CC}{\mbox{\bf C}}

\newcommand{\ZZ}{\mbox{\bf Z}}

\begin{document}

\maketitle

\begin{abstract} 
The paper aims to generalize Clausen's identity to the square of any Gauss hypergeometric function.
Accordingly, solutions of the related 3rd order linear differential equation are found in terms of
certain bivariate series that can reduce to $\hpgo32$ series similar to those in Clausen's identity.
The general contiguous variation of Clausen's identity is found. The related Chaundy's identity 
is generalized without any restriction on the parameters of Gauss hypergeometric function.
The special case of dihedral Gauss hypergeometric functions is underscored.
\end{abstract}

\section{Introduction}
\label{sec:symsq}

Clausen's identity \cite{clausen28} is
\begin{equation} \label{eq:clausen}
\hpg21{a,\;b}{a+b+\frac12}{\,z}^2=\hpg32{2a,\,2b,\,a+b}{2a+2b,a+b+\frac12}{\,z}.
\end{equation}
A related identity is \cite[(5)]{chaundy58}: 
\begin{equation} \label{eq:altclaus}
\hpg21{a,\;b}{a+b+\frac12}{\,z}\hpg21{\frac12-a,\frac12-b}{\frac32-a-b}{\,z}=
\hpg32{\frac12,a-b+\frac12,b-a+\frac12}{a+b+\frac12,\,\frac32-a-b}{\,z}.
\end{equation}
Up to a power factor, the two $\hpgo21$ functions in the second identity are solutions
of the same Euler's hypergeometric differential equation. The two identities demonstrate a case
when the {\em symmetric tensor square} of Euler's hypergeometric equation coincides
with the hypergeometric differential equation for $\hpgo32$ functions. 
Accordingly, formulas (\ref{eq:clausen}) and (\ref{eq:altclaus})
relate quadratic forms in $\hpgo21$ functions and $\hpgo32$ functions as solutions
of the same third order linear differential equation; see Exercise 13 in \cite[pg.~116]{specfaar}.

This paper aims to generalize Clausen's identity to the square of a general $\hpgo21$ function,
without any restriction on the three parameters. Accordingly, we try to find other attractive solutions
of the symmetric tensor square equation for general Euler's hypergeometric equation.
We present the symmetric tensor square equation in Section \ref{sec:deqs}; 
see formula (\ref{eq:symsq}).

We find that the general symmetric square equation has solutions expressible
as specializations of the following double hypergeometric series:
\begin{eqnarray} \label{eq:gap211}
F^{2:1;1}_{1:1;1}\!\left(\left. {a; b; p_1, p_2\atop c; q_1, q_2} \right| x,y\right) \equal
\sum_{i=0}^{\infty} \sum_{j=0}^{\infty}
\frac{(a)_{i+j}\,(b)_{i+j}\,(p_1)_i\,(p_2)_j}{(c)_{i+j}\,(q_1)_i\,(q_2)_j\;i!\,j!}\,x^i\,y^j,\\
\label{eq:gap122}
F^{1:2;2}_{2:0;0}\!\left(\left. {a; p_1, p_2; q_1, q_2 \atop b; \; c} \right| x,y\right) \equal
\sum_{i=0}^{\infty} \sum_{j=0}^{\infty}
\frac{(a)_{i+j}\,(p_1)_i\,(p_2)_j\,(q_1)_i\,(q_2)_j}{(b)_{i+j}\;(c)_{i+j}\;i!\,j!}\,x^i\,y^j.
\end{eqnarray}
If $a=c$, these two series reduce to the bivariate Appell's $F_2$ or $F_3$ hypergeometric series, respectively.
 
Surely, the left-hand sides of (\ref{eq:clausen})--(\ref{eq:altclaus}) are trivially double hypergeometric series as well. The significance of new solutions is that they reduce to $\hpgo32$ functions
of similar shape as in (\ref{eq:clausen})--(\ref{eq:altclaus}) in special cases, 
and that they represent elementary solutions (that is, a polynomial times a power function)
if the the symmetric square equation has such solutions.

Our main results are summarized in the following two theorems. In Theorem \ref{th:sameeq},
by a univariate specialization of $F^{2:1;1}_{1:1;1}$ or $F^{2:1;1}_{1:1;1}$ function we understand
the restriction of any branch (under analytic continuation) of the bivariate 
$F^{2:1;1}_{1:1;1}$ or $F^{2:1;1}_{1:1;1}$ function to the curve parameterized by $z$. 
In Theorem \ref{th:gclausen}, we refer to the $F^{2:1;1}_{1:1;1}$ or $F^{2:1;1}_{1:1;1}$ series
as defined in (\ref{eq:gap211})--(\ref{eq:gap122}). 
A direct generalization of Clausen's identity is formula (\ref{ge:clausen}). 
It applies to the Gauss hypergeometric functions contiguous to Clausen's instance.
We prove these theorems in Sections \ref{sec:sameeq} and \ref{sec:gclausen}, respectively.
\begin{theorem} \label{th:sameeq}
The univariate functions 
\begin{eqnarray}  \label{eq:spap0} 
F^{2:1;1}_{1:1;1}\!\left(\left. {2a;\,2b;\;c-\frac12,\,a+b-c+\frac12 
\atop a+b+\frac12;\, 2c-1, 2a+2b-2c+1} \right|  z,1-z\right),\\
\label{eq:spap1} z^{1-c}(1\!-\!z)^{c-a-b-\frac12}
F^{1:2;2}_{2:0;0}\!\left(\left. {\!\frac12; a\!-\!b\!+\!\frac12,a\!+\!b\!-\!c\!+\!\frac12; 
b\!-\!a\!+\!\frac12, c\!-\!a\!-\!b\!+\!\frac12 \atop c;\,2-c} 
\right|  z,\frac{z}{z\!-\!1}\right),\\ \label{eq:spap2}
z^{\frac12-c}(1\!-\!z)^{c-a-b}\,
F^{1:2;2}_{2:0;0}\!\left(\left. {\!\frac12;\, a-b+\frac12,c-\frac12;\, 
b-a+\!\frac12, \frac32-c \atop a+b-c+1;\,c-a-b+1} 
\right|  1-z,1-\frac{1}{z}\right),\\ \label{eq:spap3}
z^{\frac12-c}(1\!-\!z)^{c-a-b-\frac12}
F^{1:2;2}_{2:0;0}\!\left(\left. {\!\frac12; c-\frac12,a+b-c+\frac12; 
\frac32-c, c-a-b+\frac12 \atop a-b+1;\,b-a+1} 
\right|  \frac1z,\frac{1}{1\!-\!z}\right)
\end{eqnarray}
satisfy the symmetric tensor square equation for $\hpg21{\!a,\,b}c{z}^2$.
\end{theorem}

\begin{theorem} \label{th:gclausen}
The following identities hold in a neighborhood of $z=0$:
\begin{eqnarray} \label{ge:altclaus}
\rhallign  \hpg{2}{1}{a,\,b\,}{c}{\,z\,}\,\hpg{2}{1}{1+a-c,\,1+b-c}{2-c}{\,z\,}=
(1-z)^{c-a-b-\frac12}\times \nonumber\\
\rhallignh F^{1:2;2}_{2:0;0}\!\left(\left. {\frac12; a-b+\frac12,a+b-c+\frac12; 
b-a+\frac12, c-a-b+\frac12 \atop c;\;2-c} 
\right|  z,\frac{z}{z\!-\!1}\right),\\ \label{ge:clausen}
\rhallign \hpg21{a,\;b}{a+b+n+\frac12}{\,z}^2=  \nonumber\\ 
\rhallignh \frac{(\frac12)_n\,(a+b+\frac12)_n}{(a+\frac12)_n\,(b+\frac12)_n}\,
F^{2:1;1}_{1:1;1}\!\left(\left. {2a;\,2b;\;a+b+n,\,-n 
\atop a+b+\frac12;\, 2a+2b+2n, -2n} \right|  z,1-z\right),\\
\rhallign \frac{\left(a+\frac12\right)_n\left(a+m+\frac12\right)_n}
{\left(\frac12\right)_n\left(m+\frac12\right)_n}\,
\hpg21{a,\;-a-m-n}{\frac12-m}{\,z}^2+ \nonumber\\  \label{ge:dclausen}
\rhallign \frac{\left(a+\frac12\right)_{m}\!\left(a+n+\frac12\right)_m(a)_{m+n+1}^2}
{\left(\frac12\right)_m\left(\frac12\right)_{m+1}^2\left(\frac12\right)_n\left(\frac12\right)_{m+n}}
\,z^{2m+1}\,\hpg21{a+m+\frac12,\frac12-a-n}{m+\frac32}{\,z}^2 \nonumber\\
\rhallignh =F^{2:1;1}_{1:1;1}\!\left(\left. {2a;\,-2a-2m-2n;\;-m,\,-n 
\atop \frac12-m-n;\, -2m, -2n} \right|  z,1-z\right).
\end{eqnarray}
Here $a,b,c$ can be any complex numbers, if only the involved lower parameters $c$, $2-c$,
$a+b+\frac12$ are not zero or negative integers.
The numbers $m,n$ are assumed to be non-negative integers. 
\end{theorem}
The $F^{2:1;1}_{1:1;1}$ series in $(\ref{ge:clausen})$ is understood to terminate
in the second argument at the power $(1-z)^n$, but it is never a terminating series in
the first argument even if $a+b+n$ is a non-positive integer $-m$ (when the term-wise limit
$b\to -a-n-m$ should be applied). The $F^{2:1;1}_{1:1;1}$ series 
in $(\ref{ge:dclausen})$ is a finite rectangular sum with $(m+1)(n+1)$ terms.

If $n=0$, formula (\ref{ge:clausen}) reduces to Clausen's identity (\ref{eq:clausen}).
If $c=a+b+\frac12$, formula (\ref{ge:altclaus}) reduces to (\ref{eq:altclaus}). 
 We should observe that Clausen's identity is wrong if $a+b$ is a non-positive integer
 and the $\hpgo32$ series is interpreted as terminating. A correct identity with the
 terminating $\hpgo32$ series is then the special case $n=0$ of (\ref{ge:dclausen}):
 \begin{eqnarray} \label{eq:exclausen}
\rhallign \hpg21{a,\;-a-m}{\frac12-m}{\,z}^2+
\frac{\left(a+\frac12\right)_{m}^2(a)_{m+1}^2}
{\left(\frac12\right)_m^2\left(\frac12\right)_{m+1}^2}
\,z^{2m+1}\,\hpg21{a+m+\frac12,\frac12-a}{m+\frac32}{\,z}^2 \nonumber\\
\rhallignh =\hpg32{2a,\,-2a-2m,\,-m}{\frac12-m,\, -2m}{\,z},
\end{eqnarray}

Formula (\ref{ge:dclausen}) has the following significance. The $\hpgo21$ functions 
on the left-hand side are {\em dihedral Gauss hypergeometric functions},
as the monodromy group of their Euler's hypergeometric equation is a dihedral group.
Then the symmetric tensor square equation is reducible, and its monodromy representation
has an invariant one-dimensional subspace. 
The $F^{2:1;1}_{1:1;1}$  function on the right-hand side is a polynomial, so it is obviously
an invariant of the monodromy group. Formula (\ref{ge:dclausen}) identifies a generator for
the invariant space as a linear combination of $\hpgo21(z)^2$ solutions, and as an explicit
terminating $F^{2:1;1}_{1:1;1}$ expression \cite[Section 5]{tdihedral}.
In more plain terms, the $F^{2:1;1}_{1:1;1}$ expression gives an expected elementary solution
of the symmetric square equation. Similarly, the $F^{1:2;2}_{2:0;0}$ series in (\ref{ge:altclaus}) 
is terminating in both summation directions in the following specialization
with dihedral $\hpgo21$ functions:
\begin{eqnarray} \label{ge:altclaust}
\hpg{2}{1}{a,\,a+m+\frac12}{2a+m+n+1}{\,z\,}\hpg{2}{1}{-a-m-n,\,\frac12-a-n}{1-2a-m-n}{\,z\,}=
(1-z)^{n}\times\hspace{28pt}  \nonumber\\
F^{1:2;2}_{2:0;0}\!\left(\left. {\frac12;\; -m,-n;\,
m+1, n+1 \atop 2a+m+n+1;1-2a-m-n} 
\right|  z,\frac{z}{z\!-\!1}\right).
\end{eqnarray}
The terminating $F^{2:1;1}_{1:1;1}$ and $F^{1:2;2}_{2:0;0}$ sums are closely related.
Up to a constant multiple, one can rewrite the terminating $F^{2:1;1}_{1:1;1}$ sum 
as a terminating $F^{1:2;2}_{2:0;0}$ sum,
and vice versa, simply by reversing the order of summation in both directions.
 

\section{Variation of formulas}

A product of two single hypergeometric series is trivially a double hypergeometric series. 
In particular, we can write
\begin{equation} \label{eq:hpgabc2}
\hpg21{a,\;b}{c}{\,z}^2= 
F^{0:2;2}_{0:1;1}\!\left(\left. {a,\,a;\,b,\,b \atop c,\;c} \right|  z,\,z\right).
\end{equation}
Chaundy \cite[(6)]{chaundy58} gives the following double series expansion:
\begin{eqnarray} \label{eq:gclau2}
\hpg21{a,\;b}{c}{\,z}^2=\sum_{k=0}^{\infty}
\frac{(2a)_k(2b)_k(c\!-\!\frac12)_k}{(c)_k\,(2c-1)_k\,k!}
\hpg43{\!-\frac{k}{2},-\frac{k-1}2,\frac12,a\!+\!b\!-\!c\!+\!\frac12}{a+\frac12,\,b+\frac12,\,\frac32-c-k}{\,1}z^k.
\end{eqnarray}
But our formulas (\ref{ge:clausen}) and (\ref{ge:altclaus}) are particularly interesting as 
direct generalizations of the classical formulas (\ref{eq:clausen})--(\ref{eq:altclaus}).

Application of Euler's transformation \cite[(2.2.7)]{specfaar} to the second $\hpgo21$ factor
on the left-hand side of  (\ref{ge:altclaus}) gives
\begin{eqnarray}  \label{ge:altclaus2}
\rhallign \hpg{2}{1}{a,\,b\,}{c}{\,z\,}\,\hpg{2}{1}{1-a,\,1-b}{2-c}{\,z\,}=
\frac1{\sqrt{1-z}}\times \nonumber\\ 
\rhallignh F^{1:2;2}_{2:0;0}\!\left(\left. {\!\frac12; a-b+\frac12,a+b-c+\frac12; 
b-a+\frac12, c-a-b+\frac12 \atop c;\,2-c} 
\right|  z,\frac{z}{z\!-\!1}\right).
\end{eqnarray}
The substitution $z\mapsto z/(z-1)$ and application of Pfaff's transformation \cite[(2.2.6)]{specfaar}
to the $\hpgo21$ functions in (\ref{ge:dclausen}) and (\ref{eq:exclausen}) gives, respectively:
\begin{eqnarray} \label{ge:dclausen2}
\rhallign \frac{\left(a+\frac12\right)_n\left(a+m+\frac12\right)_n}
{\left(\frac12\right)_n\left(m+\frac12\right)_n}\,
\hpg21{a,\;a+n+\frac12}{\frac12-m}{\,z}^2- \nonumber\\  
\rhallign \frac{\left(a+\frac12\right)_{m}\!\left(a+n+\frac12\right)_m(a)_{m+n+1}^2}
{\left(\frac12\right)_m\left(\frac12\right)_{m+1}^2\left(\frac12\right)_n\left(\frac12\right)_{m+n}}
\,z^{2m+1}\,\hpg21{a+m+\frac12,a+n+m+1}{m+\frac32}{\,z}^2 \nonumber\\
\rhallignh =(1-z)^{-2a}\,F^{2:1;1}_{1:1;1}\!\left(\left. {2a;\,-2a-2m-2n;\;-m,\,-n 
\atop \frac12-m-n;\, -2m, -2n} \right|  \frac{z}{z-1},\frac1{1-z}\right),\\
\rhallign \hpg21{a,\;a+\frac12}{\frac12-m}{\,z}^2-
\frac{\left(a+\frac12\right)_{m}^2(a)_{m+1}^2}
{\left(\frac12\right)_m^2\left(\frac12\right)_{m+1}^2}
\,z^{2m+1}\,\hpg21{a+m+\frac12,a+m+1}{m+\frac32}{\,z}^2 \nonumber\\
\rhallignh =(1-z)^{-2a}\,\hpg32{2a,\,-2a-2m,\,-m}{\frac12-m,\, -2m}{\frac{z}{z-1}}.
\end{eqnarray}
Notice that
\begin{equation}
\frac{\left(a+\frac12\right)_{m}(a)_{m+1}}{\left(\frac12\right)_m\left(\frac12\right)_{m+1}}
=\frac{2^{2m}\,(m!)^2\,(2a)_{2m+1}}{(2m)!\,(2m+1)!}.
\end{equation}
The specializations $m=0$ of (\ref{ge:dclausen}) and (\ref{ge:dclausen2}) are interesting as well:
\begin{eqnarray}
\rhallign {\textstyle\left(a+\frac12\right)_n^2} \,\hpg21{a,\;-a-n}{1/2}{\,z}^2
+2\,z\,(a)_{n+1}^2\,\hpg21{a+\frac12,\frac12-a-n}{3/2}{\,z}^2 \nonumber\\
\rhallignh = {\textstyle \left(\frac12\right)^2_n}\;
\hpg32{2a,\,-2a-2n,\,-n}{\frac12-n,\,-2n}{1-z},\\
\rhallign {\textstyle\left(a+\frac12\right)_n^2} \,\hpg21{a,\;a+n+\frac12}{1/2}{\,z}^2
-2\,z\,(a)_{n+1}^2\,\hpg21{a+\frac12,a+n+1}{3/2}{\,z}^2 \nonumber\\
\rhallignh = {\textstyle \left(\frac12\right)^2_n}\;(1-z)^{-2a}\,
\hpg32{2a,\,-2a-2n,\,-n}{\frac12-n,\,-2n}{\frac1{1-z}}.
\end{eqnarray}

It is tempting to equate $\hpg21{\!a,\,b}c{z}^2$ to the 
$F^{2:1;1}_{1:1;1}$ function in (\ref{eq:spap0}) up to a constant multiple,
since the latter appears to have the power series expansion at $z=0$ as well. 
The $k$th derivative of the $F^{2:1;1}_{1:1;1}$ function at $z=0$ would evaluate
to a linear combination of the values
\begin{equation}
\hpg32{2a+k,\,2b+k,\,a+b-c+j+\frac12}{a+b+k+\frac12,\,2a+2b-2c+j+1}{\,1\,}, \qquad j=0,1,\ldots,k.
\end{equation}
However, the convergence condition of these $\hpgo32(1)$ sums is $\mbox{Re}(1-c-k)>0$,
so the $\hpgo32(1)$ values are undefined for large enough $k$. 
The $F^{2:1;1}_{1:1;1}$ function 
can have branching behavior at $z=0$ in general, as the local exponents 
$z^{1-c}$, $z^{2-2c}$ of the symmetric square equation can come into play.
See \cite{bisspec} for explicit details.

The $F^{2:1;1}_{1:1;1}$ function in (\ref{eq:spap0}) can be evaluated at $z=0$ and $z=1$
using the $\hpgo32(1)$ evaluation in \cite[Theorem 3.5.5{\em (i)}]{specfaar}:
\begin{eqnarray}
F^{2:1;1}_{1:1;1}\!\left(\left. {2a;\,2b;\;c-\frac12,\,a+b-c+\frac12 
\atop a+b+\frac12;\, 2c-1, 2a+2b-2c+1} \right|  0,1\right)\equal \nonumber\\
&&\hspace{-110pt}
\frac{\Gamma\left(\frac12\right)\Gamma\left(a+b+\frac12\right)\Gamma(1-c)\,\Gamma(1+a+b-c)}
{\Gamma\left(a+\frac12\right)\Gamma\left(b+\frac12\right)\Gamma(1+a-c)\,\Gamma(1+b-c)},\\
F^{2:1;1}_{1:1;1}\!\left(\left. {2a;\,2b;\;c-\frac12,\,a+b-c+\frac12 
\atop a+b+\frac12;\, 2c-1, 2a+2b-2c+1} \right|  1,0\right)\equal \nonumber\\
\label{eq:gap210}&&\hspace{-110pt}
\frac{\Gamma\left(\frac12\right)\Gamma\left(a+b+\frac12\right)\Gamma(c)\,\Gamma(c-a-b)}
{\Gamma\left(a+\frac12\right)\Gamma\left(b+\frac12\right)\Gamma(c-a)\,\Gamma(c-b)}.
\end{eqnarray}
The convergence conditions are $\mbox{Re}(1-c)>0$ and $\mbox{Re}(c-a-b)>0$,
respectively.

It is worth mentioning that Bailey's identity \cite{Bailey33} 
\begin{equation} \label{eq:bailey}
\app4{a;\;b}{c,a\!+\!b\!-\!c\!+\!1}{x(1-y),y(1-x)}
=\hpg21{a,b}{c}{x} \hpg21{a,\;b}{a\!+\!b\!-\!c\!+\!1}{y}
\end{equation}
specializes to 
\begin{equation} \label{eq:bailey2}
\app4{a;\;b}{c,a\!+\!b\!-\!c\!+\!1}{x^2,(1-x)^2}
=\hpg21{a,b}{c}{x} \hpg21{a,\;b}{a\!+\!b\!-\!c\!+\!1}{1-x}.
\end{equation}
The $F_4(x^2,(1-x)^2)$ function is another solution of
the symmetric tensor square equation. 
In particular, we have the following  identity for $x\in[0,1]$ 
if $\mbox{Re} c<1$ and $\mbox{Re}(c-a-b)>0$, 
due to the connection formula in \cite[(2.3.13)]{specfaar}:
\begin{eqnarray}
\hspace{-15pt}\hpg21{a,\,b}{c}{z}^2\!\equal
\frac{\Gamma(c)\Gamma(c-a-b)}{\Gamma(c-a)\Gamma(c-b)}\,
\app4{a;\;b}{\!c,a\!+\!b\!-\!c\!+\!1}{x^2,(1-x)^2}+ \nonumber\\
&&\!\frac{\Gamma(c)\Gamma(a\!+\!b\!-\!c)}{\Gamma(a)\,\Gamma(b)}\,(1-t)^{2c-2a-2b}
\app4{c-a;\;c-b}{\!c,c\!-\!a\!-\!b\!+\!1}{x^2,(1-x)^2}.
\end{eqnarray}

\section{The differential equations}
\label{sec:deqs}

Let $\ediff{z}$, $\ediff{x}$, $\ediff{y}$ denote the differential operators
\begin{equation} \label{eq:ediff}
\ediff{z}=z\,\frac{d}{dz}, \qquad \ediff{x}=x\,\frac{\partial}{\partial x}, \qquad
\ediff{y}=y\,\frac{\partial}{\partial y}.
\end{equation}
Euler's hypergeometric differential equation \cite[(2.3.5)]{specfaar} for a 
general Gauss hypergeometric function $\displaystyle\hpg{2}{1}{a,\,b\,}{c}{\,z\,}$
can be compactly written as follows: 
\begin{equation} \label{hpgde}
z\,(\ediff{z}+a)\,(\ediff{z}+b)\,y(z)-\ediff{z}\,(\ediff{z}+c-1)\,y(z)=0.
\end{equation}
This is a Fuchsian equation 
with three regular singular points $z=0,1$ and $\infty$. 
The local exponents at the singular points are:
\[
\mbox{$0$, $1-c$ at $z=0$;} \qquad \mbox{$0$, $c-a-b$ at $z=1$;}\qquad
\mbox{and $a$, $b$ at $z=\infty$.}
\]
In general, a basis of solutions is 
\begin{equation}
\hpg{2}{1}{a,\,b\,}{c}{\,z\,}, \qquad  z^{1-c}\,\hpg{2}{1}{1+a-c,\,1+b-c}{2-c}{\,z\,}.
\end{equation}
Note that the left-hand sides of formulas (\ref{eq:altclaus}) and (\ref{ge:altclaus}) 
are products of these two solutions of the same Euler's equation,
with the power factor $z^{1-c}$ ignored. The classical formulas  (\ref{eq:clausen})--(\ref{eq:altclaus}) apply  when the difference of local exponents at $z=1$ is equal to $1/2$.

If $y_1,y_2$ are two independent solutions of a second order homogeneous linear 
differential equation like (\ref{hpgde}), the functions $y_1^2$, $y_1y_2$, $y_2^2$ 
satisfy a third order linear differential equation called the {\em symmetric tensor square}
of the second order equation. The symmetric tensor square equation for Euler's equation
(\ref{hpgde}) can be written as follows \cite[(2)]{chaundy58}: 
\begin{eqnarray} \label{eq:symsq}
z\,(\ediff{z}+2a)\,(\ediff{z}+2b)\,(\ediff{z}+a+b)\,y(z)
-\ediff{z}\,(\ediff{z}+c-1) (\ediff{z}+2c-2) \,y(z) \nonumber\\
+\frac{(2a+2b-2c+1)\,z}{z-1}\,\big((a+b-c+1)\,\ediff{z}+2ab\big)\,y(z) \equal 0.\qquad
\end{eqnarray}
If $c=a+b+\frac12$, this is a differential equation for the $\hpgo32$ function
in (\ref{eq:clausen}).

Partial differential equations for the 
$F^{2:1;1}_{1:1;1}(x,y)$ function in (\ref{eq:gap211})
can be obtained by considering the first order recurrence relations between its coefficients.
Writing the $F^{2:1;1}_{1:1;1}$ sum as $\sum_{i=0}^\infty\sum_{j=0}^\infty c_{i,j}$, 
we have the relations
\begin{eqnarray*}
\frac{c_{i+1,j}}{c_{i,j}}=\frac{(a+i+j)(b+i+j)(p_1+i)\,x}{(c+i+j)\,(q_1+i)\,(1+i)}, \quad
\frac{c_{i,j+1}}{c_{i,j}}=\frac{(a+i+j)(b+i+j)(p_2+j)\,y}{(c+i+j)\,(q_2+i)\,(1+j)},
\end{eqnarray*}
also $c_{i+1,j}\big/c_{i,j+1}=(p_1+i)(q_2+j)(1+j)\,x\big/(p_2+j)(q_1+i)(1+i)\,y$.
These relations translate to the following partial differential equations for 
the $F^{2:1;1}_{1:1;1}(x,y)$ function:
\begin{eqnarray} \label{eq:pdo1}
P_1\equal x\,(\ediff{x}+\ediff{y}+a)\,(\ediff{x}+\ediff{y}+b)\,(\ediff{x}+p_1) \nonumber \\
&& \qquad\qquad -\ediff{x}\,(\ediff{x}+\ediff{y}+c-1)\,(\ediff{x}+q_1-1), \\  \label{eq:pdo2}
P_2\equal y\,(\ediff{x}+\ediff{y}+a)\,(\ediff{x}+\ediff{y}+b)\,(\ediff{y}+p_2) \nonumber \\
&& \qquad\qquad -\ediff{y}\,(\ediff{x}+\ediff{y}+c-1)\,(\ediff{y}+q_2-1), \\  \label{eq:pdo3}
P_3\equal x\,\ediff{y}\,(\ediff{y}+q_2-1)\,(\ediff{x}+p_1)
-y\,\ediff{x}\,(\ediff{x}+q_1-1)\,(\ediff{y}+p_2).
\end{eqnarray}
In general, two of these operators 
generate the ideal $J$ in $\CC(x,y)\langle\ediff{x},\ediff{y}\rangle$ annihilating the
$F^{2:1;1}_{1:1;1}(x,y)$ series. In particular, we have the following obvious syzygies:
\begin{eqnarray*}
-y\,(\ediff{y}+p_2)\,P_1+x\,(\ediff{x}+p_1)\,P_2
+{\textstyle (\ediff{x}+\ediff{y}+c-2)\,P_3} \equal 0,\\
\ediff{y}(\ediff{y}+q_2\!-\!1)P_1-\ediff{x}(\ediff{x}+q_1\!-\!1)P_2 
+(\ediff{x}\!+\ediff{y}\!+a-\!1)(\ediff{x}\!+\ediff{y}\!+b-\!1)P_3 \equal 0.
\end{eqnarray*}
If the coefficient $\ediff{x}+\ediff{y}+c-2$ does not divide 
$(\ediff{x}+\ediff{y}+a-1)(\ediff{x}+\ediff{y}+b-1)$, we can express
$P_3$ in terms of $P_1,P_2$, etc. 


Note the following commutation relations in $\CC(x,y)\langle\ediff{x},\ediff{y}\rangle$:
\begin{equation}
\ediff{x}x=x\ediff{x}+x=x\,(\ediff{x}+1),\qquad
\ediff{y}y=y\,(\ediff{y}+1).
\end{equation}
The variables $x,\ediff{x}$ commute with $y,\ediff{y}$. Up to the multiplication order, 
these relations are the same as shift operator relations, such as $S_nn=(n+1)S_n$
in the Ore algebra $\CC\langle n,S_n\rangle$. 
We will keep working with Euler-type differentiation operators in (\ref{eq:ediff})
in our transformations of differential equations.

Gr\"obner basis computations show that third order partial differential operators in the ideal $J$
are linearly generated by $P_1,P_2,P_3$. A residue basis over $\CC(x,y)$ is formed by 
the 7 monomials $1,\ediff{x},\ediff{y},\ediff{x}^2,\ediff{x}\ediff{y},\ediff{y}^2,\ediff{x}^3$.
Hence the rank of the differential system for $F^{2:1;1}_{1:1;1}(x,y)$ functions is equal to $7$.
The leading coefficients in various Gr\"obner bases suggests that the following lines
are in the singular locus of the differential system:
\begin{eqnarray*}
x=0,\quad x=1,\quad y=0,\quad y=1,\quad x+y=0,\quad x+y=1,\quad
y+2x=1,\quad\,2y+x=1.
\end{eqnarray*}

It is easy to apply the transformations $F(x,y)\mapsto x^\alpha y^\beta F(x,y)$ or
$x\mapsto 1/x$, $y\mapsto 1/y$ to the differential system 
generated by the operators $P_1,P_2$. In effect, we only have to replace additionally
$\ediff{x}\mapsto\ediff{x}+\alpha$, $\ediff{y}\mapsto\ediff{y}+\beta$ or
$\ediff{x}\mapsto-\ediff{x}$, $\ediff{y}\mapsto-\ediff{y}$, respectively.
We may attempt to find transformations of the operators $P_ 1,P_2$ to a pair of similar 
hypergeometric operators. This gives us the following solutions of the same system 
of partial differential equations, beside the $F^{2:1;1}_{1:1;1}(x,y)$ function in (\ref{eq:gap211}):
\begin{eqnarray}
x^{1-q_1}\,F^{2:1;1}_{1:1;1}\!\left(\left. {1+a-q_1; 1+b-q_1; 1+p_1-q_1, p_2\atop 
1+c-q_1;\, 2-q_1, \, q_2} \right| x,y\right), \\
y^{1-q_2}\,F^{2:1;1}_{1:1;1}\!\left(\left. {1+a-q_2;\,1+b-q_2; p_1, 1+p_2-q_2\atop 
1+c-q_2; \,q_1,\, 2-q_2} \right| x,y\right), \\
x^{1-q_1}\,y^{1-q_2}\,F^{2:1;1}_{1:1;1}\!\!\left(\left. {2\!+\!a\!-\!q_1\!-\!q_2; 2\!+\!b\!-\!q_1\!-\!q_2; 
1\!+\!p_1\!-\!q_1, 1\!+\!p_2\!-\!q_2 \atop 2+c-q_1-q_2;\, 2-q_1, 2-q_2} \right| x,y\right), \\
x^{-p_1}\,y^{-p_2}\,F^{1:2;2}_{2:0;0}\!\left(\left. {1+p_1+p_2-c; 1+p_1-q_1, 1+p_2-q_2; 
p_1, p_2 \atop 1+p_1+p_2-a; \; 1+p_1+p_2-b} \right| \frac1x,\frac1y\right). 
\end{eqnarray}
As we see, the systems of partial differential equations for the $F^{2:1;1}_{1:1;1}$ and 
$F^{1:2;2}_{2:0;0}$ functions are easily transformable to each other. 
In the other direction, here is a $F^{2:1;1}_{1:1;1}$ function satisfying the same partial
differential equations as the $F^{1:2;2}_{2:0;0}(x,y)$ function in (\ref{eq:gap122}):
\begin{eqnarray} \label{eq:gap2112}
x^{-p_1}\,y^{-p_2}\,F^{2:1;1}_{1:1;1}\!\left(\left. {1+p_1+p_2-a; \; 1+p_1+p_2-b;\; p_1, p_2\atop 
1+p_1+p_2-c;\; 1+p_1-q_1, 1+p_2-q_2} \right|  \frac1x,\,\frac1y\right). 
\end{eqnarray}
From here, the other $F^{2:1;1}_{1:1;1}$ companion solutions for the 
$F^{1:2;2}_{2:0;0}(x,y)$ function in (\ref{eq:gap122}) can be obtained using the symmetries
$p_1\leftrightarrow q_1$ and $p_2\leftrightarrow q_2$.

\section{Proof of Theorem \ref{th:sameeq}}
\label{sec:sameeq}
 
First we prove that the $F^{2:1;1}_{1:1;1}$ function in (\ref{eq:spap0})
satisfies the symmetric square equation (\ref{eq:symsq}) for  $\hpg21{\!a,\,b}c{z}^2$.
The ideal of partial differential operators annihilating
\begin{equation}
F^{2:1;1}_{1:1;1}\!\left(\left. {2a;\,2b;\;c-\frac12,\,a+b-c+\frac12 
\atop a+b+\frac12;\, 2c-1, 2a+2b-2c+1} \right|  x,y \right)
\end{equation}
is generated by the polynomials $P_1$, $P_2$, $P_3$ in (\ref{eq:pdo1})--(\ref{eq:pdo3})
with the parameters specialized to linear functions in $a$, $b$, $c$.
Here are the specialized differential equations:
\begin{eqnarray} \label{eq:sap1}
0 \equal 
x\,(\ediff{x}+\ediff{y}+2a)\,(\ediff{x}+\ediff{y}+2b)\,(\ediff{x}+c-\ahalf) \nonumber\\
&& \qquad\qquad -\ediff{x}\,(\ediff{x}+\ediff{y}+a+b-\ahalf)\,(\ediff{x}+2c-2), \\ \label{eq:sap2}
0 \equal 
y\,(\ediff{x}+\ediff{y}+2a)\,(\ediff{x}+\ediff{y}+2b)\,(\ediff{y}+a+b-c+\ahalf) \nonumber\\
&& \qquad\qquad -\ediff{y}\,(\ediff{x}+\ediff{y}+a+b-\ahalf)\,(\ediff{y}+2a+2b-2c), \\ \label{eq:sap3}
0 \equal 
x\,\ediff{y}\,(\ediff{y}+2a+2b-2c)\,(\ediff{x}+c-\ahalf) \nonumber\\ 
&& \qquad\qquad -y\,\ediff{x}\,(\ediff{x}+2c-2)\,(\ediff{y}+a+b-c+\ahalf).
\end{eqnarray}
The variables and derivatives are related as follows under the 
univariate specialization under consideration:
\begin{equation} \label{eq:xyspec}
x=z,\quad y=1-z,\qquad
\frac{d}{dz}=\frac{\partial}{\partial x}-\frac{\partial}{\partial y}, \qquad
\ediff{z}=\ediff{x}+\frac{z}{z-1}\,\ediff{y}.
\end{equation}
Following \cite[Definition 1.1]{AppelGhpg}, the {\em partial differential form} of an ordinary differential equation under a specialization like (\ref{eq:xyspec}) is the expression where the univariate derivatives are replaced by respective linear combinations of partial derivatives. On the other hand, 
the {\em specialized form} of a partial differential equation under the same kind of 
specialization is the expression with coefficients to the partial derivatives specialized to univariate functions.  In our setting, the two forms are ``mixed" linear differential expressions in 
the partial derivatives $\Theta_x$, $\Theta_y$ but with the coefficients univariate in $z$.
Algebraically, we work in the $\CC[z]\langle\Theta_z\rangle$ module generated by the
partial derivatives of any order. The action of $\Theta_z$ on the partial derivatives is given by
the identification in (\ref{eq:xyspec}). 

To show that the univariate $F^{2:1;1}_{1:1;1}$ function in (\ref{ge:clausen})
satisfies the symmetric tensor square equation (\ref{eq:symsq}), we demonstrate
that  the partial differential form of the symmetric square equation coincides
with the specialized form of a partial differential equation following from 
the partial differential equations (\ref{eq:sap1})--(\ref{eq:sap3}). 
In the terminology of \cite[Definition 1.1]{AppelGhpg}, we show that
the symmetric square equation {\em follows fully} from the differential equations
(\ref{eq:sap1})--(\ref{eq:sap3}) under the specialization defined in (\ref{eq:xyspec}).

To give a partial differential form of the symmetric square equation (\ref{eq:symsq}),
we introduce a way to work with multiplicative expressions of (non-commutative) differential operators.
Suppose that $G$, $H$ are functions in $z$. Then
\begin{eqnarray} \label{eq:newprodgh}
&& \hspace{-40pt} (\ediff{z}+G)\,(\ediff{z}+H) =
(\ediff{z}+G)\left(\ediff{x}+\frac{z}{z-1}\,\ediff{y}+H\right) \nonumber \\ 
&&  = \left(\ediff{x}+\frac{z}{z-1}\,\ediff{y}+A\right) \comm
\left(\ediff{x}+\frac{z}{z-1}\,\ediff{y}+B\right)-\frac{z}{(z-1)^2}\,\ediff{y}+z\,\frac{dH}{dz},
\end{eqnarray}
where we use the dot between two differential factors to signify {\em commutative} multiplication. 
In other words, an expanded expression for $(\ediff{z}+G)\,(\ediff{z}+H)$ can be obtained by
multiplying commutatively the two large factors in (\ref{eq:newprodgh}), collecting terms to
the monomials in $\Theta_x$, $\Theta_y$, then interpreting those monomials as differential operators
of suitable order, and writing the coefficients to those monomials as left-side factors. 

If $A$, $B$, $C$ are constants,  the product
$(\ediff{z}+A)\,(\ediff{z}+B)\,(\ediff{z}+C)$ can be written
as follows:
\begin{eqnarray*}
\left(\ediff{x}+\frac{z}{z-1}\ediff{y}+A\right) \comm
\left(\left(\ediff{x}+\frac{z}{z-1}\ediff{y}+B\right) \comm
\left(\ediff{x}+\frac{z}{z-1}\ediff{y}+C\right)-\frac{z}{(z-1)^2}\,\ediff{y}\right)\\
-\frac{z}{(z-1)^2} \,\ediff{y}\, \comm \left(2\ediff{x}+\frac{2z}{z-1}\ediff{y}+B+C\right)
+\frac{z\,(z+1)}{(z-1)^3}\,\ediff{y},
\end{eqnarray*}
and finally, as
\begin{eqnarray}
\left(\ediff{x}+\frac{z}{z-1}\ediff{y}+A\right) \comm
\left(\ediff{x}+\frac{z}{z-1}\ediff{y}+B\right) \comm
\left(\ediff{x}+\frac{z}{z-1}\ediff{y}+C\right) \nonumber \\
-\frac{z}{(z-1)^2} \,\ediff{y}\, \comm \left(3\ediff{x}+\frac{3z}{z-1}\ediff{y}+A+B+C\right)
+\frac{z\,(z+1)}{(z-1)^3}\,\ediff{y}.
\end{eqnarray}
Following this ``commutative" expression,
we can write the symmetric square equation (\ref{eq:symsq}) as follows:
\begin{eqnarray}
&&  \hspace{-26pt} z\,\left(\ediff{x}+\frac{z}{z-1}\ediff{y}+2a\right) \comm
\left(\ediff{x}+\frac{z}{z-1}\ediff{y}+2b\right) \comm
\left(\ediff{x}+\frac{z}{z-1}\ediff{y}+a+b\right) \nonumber \\
&& -\left(\ediff{x}+\frac{z}{z-1}\ediff{y}\right) \comm
\left(\ediff{x}+\frac{z}{z-1}\ediff{y}+c-1\right) \comm
\left(\ediff{x}+\frac{z}{z-1}\ediff{y}+2c-2\right) \hspace{-25pt} \nonumber \\
&& -\frac{3z}{(z-1)^2} \,\ediff{y}\, \comm \left(\ediff{x}+\frac{z}{z-1}\ediff{y}+\frac{(a+b)\,z+1-c}{z-1}\right)
+\frac{z\,(z+1)}{(z-1)^2}\,\ediff{y} \nonumber \\
&& +\frac{(2a+2b-2c+1)\,z}{z-1}\,\left((a+b-c+1)\left(\ediff{x}+\frac{z}{z-1}\ediff{y}\right)+2ab\right)=0.
\hspace{-18pt}
\end{eqnarray}
Using formal commutative computations, one can check that this equation
coincides with the specialized (under $x=z$, $y=1-z$) form of
\begin{equation}
\mbox{[ Eq.~(\ref{eq:sap1}) ]} -\frac{z^2}{(z-1)^2}\mbox{ [ Eq.~(\ref{eq:sap2}) ]}
-\frac{1}{(z-1)^2}\mbox{ [ Eq.~(\ref{eq:sap3}) ]}.
\end{equation}
We proved that the $F^{2:1;1}_{1:1;1}$ function in (\ref{eq:spap0})
satisfies the symmetric square equation (\ref{eq:symsq}).

The relation between the functions in (\ref{eq:gap122}) and (\ref{eq:gap2112})
implies that the function in (\ref{eq:spap3}) satisfies the same symmetric square equation.
The symmetric square equation has the following two solutions as well:
\begin{equation}
z^{-2a}\,\hpg21{a,\;1+a-c}{1+a-b}{\,\frac1z\,}^2, \qquad 
(1-z)^{-2a}\,\hpg21{a,\;c-b}{1+a-b}{\frac1{1-z}}^2.
\end{equation}
The relation between the functions $\hpg21{\!a,\,b}c{z}^2$ and (\ref{eq:spap3})
translates into the claim that (\ref{eq:spap1}) and (\ref{eq:spap2}) are solutions 
of the same symmetric square equation as well. This completes the proof of
Theorem \ref{th:sameeq}.

Notice that the upper parameters in  (\ref{eq:spap1})--(\ref{eq:spap3}) contain the local exponent differences of $\hpg21{\!a,\,b}c{z}$ increased by $1/2$, while the lower parameters are the same local exponent differences increased by $1$. A set of twelve $F^{2:1;1}_{1:1;1}$ solutions can be obtained using the relation between the functions in (\ref{eq:gap122}) and (\ref{eq:gap2112}). 

\section{Proof of Theorem \ref{th:gclausen}}
\label{sec:gclausen}

To prove formula (\ref{ge:altclaus}), we compare the function in (\ref{eq:spap1}) with
the following solution of the symmetric square equation:
\begin{equation} \label{eq:ssqsolp}
z^{1-c}\;\hpg{2}{1}{a,\,b\,}{c}{\,z\,}\,\hpg{2}{1}{1+a-c,\,1+b-c}{2-c}{\,z\,}.
\end{equation}
Both solutions have the same local exponent $1-c$ at $z=0$.
For those values of $c$ for which both $\hpgo21$ solutions are generally defined,
the linear space of solutions with this local exponent is one-dimensional.
After division by $z^{1-c}$ both solutions (\ref{eq:spap1}) and (\ref{eq:ssqsolp}) evaluate
to $1$ at $z=0$, therefore formula (\ref{ge:altclaus}) follows.

To prove formula (\ref{ge:clausen}), we observe that (for general $a,b$)
the linear space of power series solutions at $z=0$ to the symmetric square equation 
is one-dimensional. The $F^{2:1;1}_{1:1;1}$ function in (\ref{ge:clausen}) and the square
function in (\ref{eq:hpgabc2}) are both power series solutions with the local exponent $0$.
Hence they must differ by a constant multiple. To compute the constant multiple,
we evaluate both sides of $z=1$. Note that the convergence condition 
$\mbox{Re}\left(n+\frac12\right)>0$ for both sides 
is always satisfied. The left-hand side is evaluated by 
the Gauss formula \cite[Theorem 2.2.2]{specfaar}; 
the right-hand side is evaluated by (\ref{eq:gap210}) with $c=a+b+n+\frac12$. 
The Pochhammer-type factor follows.

To prove formula (\ref{ge:dclausen}), we take the term-wise limit $b\to-a-m-n$ in (\ref{ge:clausen}).
We use
\begin{equation}
\lim_{\varepsilon\to 0} \frac{(\varepsilon-k)_{2k+1}}{(2\varepsilon-2k)_{2k+1}}=
\frac{(-k)_k\,k!}{(-2k)_{2k}\cdot2}=\frac{(-1)^k\,(k!)^2}{2\cdot(2k)!}
=\frac{(-1)^k\,k!}{2^{2k+1}(\frac12)_k}.
\end{equation}
to get the limiting formula
\begin{eqnarray}  \label{ge:altern2}
\frac{\left(a+\frac12\right)_n\left(\frac12-a-m-n\right)_n}
{\left(\frac12\right)_n\left(\frac12-m-n\right)_n}\,
\hpg21{a,\;-a-m-n}{\frac12-m}{\,z}^2 \hspace{-190pt} \nonumber\\ 
\equal F^{2:1;1}_{1:1;1}\!\left(\left. {2a;\,-2a-2m-2n;\;-m,\,-n 
\atop \frac12-m-n;\, -2m, -2n} \right|  z,1-z\right)+ \nonumber\\
&&\frac{(2a)_{2m+1}\,(-2a-2m-2n)_{2m+1}}{\left(\frac12-m-n\right)_{2m+1}(2m+1)!}
\frac{(-1)^m\,m!}{2^{2m+1}(\frac12)_m}\,z^{2m+1}\times\nonumber\\
&&F^{2:1;1}_{1:1;1}\!\left(\left. {2a+2m+1;\,1-2a-2n;\; m+1,\,-n 
\atop \frac32+m-n;\, 2m+2, -2n} \right|  z,1-z\right).
\end{eqnarray}
Due to (\ref{ge:clausen}), the latter $F^{2:1;1}_{1:1;1}$ function can be written as 
\begin{eqnarray}
\frac{\left(a+m+1\right)_n\left(1-a-n\right)_n}
{\left(\frac12\right)_n\left(m-n+\frac32\right)_n}\,
\hpg21{a+m+\frac12,\frac12-a-n}{m+\frac32}{\,z}^2.
\end{eqnarray}
We replace $(2m+1)!=2^{2m+1}m!\,(\frac12)_{m+1}$ and rewrite identity (\ref{ge:altern2}) as follows:
\begin{eqnarray} 
\frac{\left(a+\frac12\right)_n\left(a+m+\frac12\right)_n}
{\left(\frac12\right)_n\left(m+\frac12\right)_n}\,
\hpg21{a,\;-a-m-n}{\frac12-m}{\,z}^2 \hspace{-190pt} \nonumber\\ 
\equal F^{2:1;1}_{1:1;1}\!\left(\left. {2a;\,-2a-2m-2n;\;-m,\,-n 
\atop \frac12-m-n;\, -2m, -2n} \right|  z,1-z\right)+ \nonumber\\
&&\frac{(-1)^{m+n+1}}{2^{4m+2}(\frac12)_m(\frac12)_{m+1}}\,
\frac{(2a)_{2m+1}\,(2a+2n)_{2m+1}}{\left(\frac12-m-n\right)_{2m+1}}
\,z^{2m+1}\times\nonumber\\
&&\frac{\left(a+m+1\right)_n\left(a\right)_n}
{\left(\frac12\right)_n\left(m-n+\frac32\right)_n}\,
\hpg21{a+m+\frac12,\frac12-a-n}{m+\frac32}{\,z}^2
\end{eqnarray}
After regrouping
\begin{eqnarray*}
\textstyle
\left(\frac12-m-n\right)_{2m+1}\left(m-n+\frac32\right)_n=\left(\frac12-m-n\right)_{2m+n+1}
=(-1)^{m+n}\left(\frac12\right)_{m+n}\left(\frac12\right)_{m+1},\\
\textstyle (2a)_{2m+1}=2^{2m+1}\,(a)_{m+1}\left(a+\frac12\right)_m, \qquad
(a)_{m+1}\,(a+m+1)_n=(a)_{m+n+1}, \qquad \mbox{etc.,}
\end{eqnarray*}
we obtain (\ref{ge:dclausen}).

Alternatively, we may find the constant multiple in (\ref{ge:clausen}) by evaluating both sides
at $z=0$. The evaluation of the right-hand side can be done by Zeilberger's algorithm. 
Formula (\ref{ge:dclausen}) can be idependently proved by observing that the
linear space of power series solutions at $z=0$ to the symmetric square equation is
two-dimensional, and evaluating at the two points $z=0$ and $z=1$.

\small

\bibliographystyle{alpha}
\bibliography{../hypergeometric}

\end{document}